\documentstyle[12pt,twoside]{amsart}

\title
{The arc space of a toric variety}
\author{Shihoko Ishii}

\newcommand{\bC}{{\Bbb C}}

\newcommand{\bZ}{{\Bbb Z}}

\newcommand{\bR}{{\Bbb R}}
\newcommand{\bN}{{\Bbb N}}
\newcommand{\bA}{{\Bbb A}}

\newcommand{\Spec}{\operatorname{Spec}}

\newcommand{\Hom}{\operatorname{Hom}}
\newcommand{\Sing}{\operatorname{Sing}}
\newcommand{\ord}{\operatorname{ord}}

\newcommand{\st}{{\Spec k[[t]]}}

\newcommand{\sTm}{{\Spec K[t]/(t^{m+1})}}
\newcommand{\sT}{{\Spec K[[t]]}}
\newcommand{\D}{{\Delta}}
\let \cedilla =\c
\renewcommand{\c}[0]{{\mathbb C}}  

\renewcommand{\o}[0]{{\mathcal O}} 

\newcommand{\isom}{{\tilde{\to}}} 
\newcommand{\image}{{\operatorname{Im}}}

\newcommand{\spec}[0]{\operatorname{Spec}}
\newcommand{\sing}[0]{\operatorname{Sing}}

\newcommand{\Cont}{\operatorname{Cont}}

\def\to {\longrightarrow}

\newtheorem{thm}{Theorem}[section]

\newtheorem{lem}[thm]{Lemma}
\newtheorem{cor}[thm]{Corollary}
\newtheorem{prop}[thm]{Proposition}
\newtheorem{problem}[thm]{Problem}

\theoremstyle{definition}
\newtheorem{defn}[thm]{Definition}

\newtheorem{say}[thm]{}
\newtheorem{exmp}[thm]{Example}

\newtheorem{rem}[thm]{Remark}

\theoremstyle{remark}

\begin{document}
\maketitle
\noindent
Department of Mathematics, Tokyo Institute of
Technology, Oh-Okayama, Meguro, Tokyo, Japan
\newline
e-mail : shihoko@@math.titech.ac.jp

\begin{abstract}
The Nash problem on arc families is affirmatively answered for a 
toric variety by Ishii and Koll\'ar's paper which also  shows the negative answer for 
general case.
The Nash problem is one of questions about the relation between arc families 
and valuations. 
In this paper,  the relation is described clearly for a toric variety.
The arc space of a toric variety 
 admits 
an action of the group scheme determined by the torus.
Each orbit on the arc space corresponds to a lattice point in the cone
and therefore corresponds to a toric valuation.
 The dominant relation among the orbits  is 
described in terms of the lattice points.
As a corollary we obtain the answer to the embedded version of the Nash problem for 
an invariant ideal on a toric variety.

\end{abstract}

\section{Introduction}

\noindent
 The concept of jet schemes and  arc space  over an algebraic variety 
 or an analytic space was introduced by Nash 
 in his preprint in 1968 
 which was later  published as \cite{nash}.
  These schemes  
  are considered as something to represent
  the nature of the singularities of the base space. 
  In fact, papers \cite{ein}, \cite{must01}, \cite{must02} 
  by  Musta\cedilla{t}\v{a}, Ein and Yasuda show that  geometric 
  properties of the jet schemes determine certain properties 
  of the singularities 
  of the base space. 
  Primarily the Nash problem posed in \cite{nash} is   based on this idea. 
  The Nash problem asks if the set of arc families through the 
  singularities corresponds bijectively to the set of the essential 
  components of resolutions of the singularities.
  Here an arc family through the singularities on \( X \)
   is a good component of \(
  \pi^{-1}(\sing X ) \) (see \ref{good} or \cite{i-k} 
  for the definition of a good component),
  where \( \pi \) is the canonical projection from the arc space to \( X \).
  The paper \cite{i-k} proves that if \( X \) is a toric variety,
  the answer to the Nash problem is ``yes'',
  while the paper also shows the negative answer for 
  general \( X \).

   In this paper, we study the structure of the  
   arc space of a toric variety defined over an algebraically closed 
   field \( k \) of arbitrary characteristic.
    We  prove that each jet scheme or arc space admits a canonical 
  action of the jet scheme or arc space of the torus.
  The arc space of a toric variety 
  becomes an almost homogeneous space by this action, 
  which means that the arc space is the 
  closure of one orbit.   
  A good component turns out to be the closure of a certain orbit
  and there is no non-good component in the arc space of a toric 
  variety.  

  Each orbit of the arc space corresponds to a lattice point of 
  the cone, therefore to a toric valuation,
   and the dominant relation of  two orbits 
  is translated to the order relation of the corresponding lattice 
  points.
  As a corollary we show the answer to the embedded version of Nash 
  problem posed by Ein, Lazarsfeld and Musta\cedilla{t}\v{a} in 
  \cite{ELM} for an invariant ideal on a toric variety.
  
  This paper is organized as follows: 
  In \S 2 we study some basic properties on jet schemes and arc spaces.
  The closed points in the arc spaces of varieties are discussed here.
  In \S 3 we introduce a stratification on  the 
  arc space of a toric variety according to the fan.
  Some basic properties of the arc space of a toric variety
  (non-existence of non-good components, irreducibility in any 
  characteristic) are proved here.                
  In \S 4 we study the orbits of the arc space of a toric variety
  by the action of the arc space of the torus.
  In \S 5 we give the answer to the embedded version of Nash problem 
  for  an invariant ideal on a toric variety.
  
  Throughout this paper the base field \( k \) is an algebraically closed 
  field of arbitrary characteristic unless otherwise stated.
  
  The author would like to thank Professors Kei-ichi Watanabe and 
  Ken-ichi Yoshida for useful discussions and 
  providing with the interesting example \ref{countable}. 
  She is also grateful to the members of the Singularity Seminar at 
  Nihon University for their stimulating discussion and the referee 
  for constructible comments. 
  The author is partly supported by Grant-In-Aid of the ministry of 
  science and education.
  


\section{Basic properties of jet schemes and the arc space}

\begin{defn}
  Let \( X \) be a scheme of finite type over \( k \)
and $K\supset k$ a field extension.
  For \( m\in \bN \) a  morphism \( \Spec K[t]/(t^{m+1})\to X \) is called an {\it \( m \)-jet}
  of \( X \) and \( \Spec K[[t]]\to X \) is called an {\it arc} of \( X \).
  We denote the closed point of \( \Spec K[[t]] \) by \( 0 \) and
  the generic point by \( \eta \).
\end{defn}

\begin{say}
\label{field}
  Let \( X \) be a scheme of finite type over \( k \).
  Let \( {\cal S}ch/k \) be the category of \( k \)-schemes  
   and \( {\cal S}et \) the category of sets.
  Define a contravariant functor  \( F_{m}: {\cal S}ch/k \to {\cal S}et \)
  by 
$$
 F_{m}(Y)=\Hom _{k}(Y\times_{\Spec k}\Spec k[t]/(t^{m+1}), X).
$$
  Then, \( F_{m} \) is representable by a scheme \( X_{m} \) of finite
  type over \( k \), that is
$$
 \Hom _{k}(Y, X_{m})\simeq\Hom _{k}(Y\times_{\Spec k}
\Spec k[t]/(t^{m+1}), X).
$$ 
   This \( X_{m} \) is called the {\it \( m \)-jet scheme} of \( X \).
  A \( K \)-valued point \( \alpha:\Spec K \to X_{m} \)   is regarded 
  as an \( m \)-jet \( \alpha:
  \Spec K[t]/(t^{m+1})
  \to X \).
  
  Let \( X_{\infty}=\varprojlim _{m}X_{m} \) and call it the 
{\it  
  arc space} of \( X \).
\( X_{\infty} \) is a scheme which is not of finite type over $k$, 
see \cite{DL}.
  Denote the canonical projection \( X_{\infty }\to X \) by \( 
  \pi \).
    A \( K \)-valued point \( \alpha:\Spec K \to X_{\infty} \)  
   is regarded as  an arc \( \alpha:\sT
  \to X \).

  Using  the representability of \( F_{m} \) 
we obtain the following universal property of $X_{\infty}$:
\end{say} 

\begin{prop}
\label{ft}
  Let \( X \) be a scheme of finite type over \( k \).
  Then
  \[ \Hom _{k}(Y, X_{\infty})\simeq\Hom _{k}(Y\widehat\times_{\Spec k}\st, X) \]
   for an arbitrary \( k \)-scheme \( Y \),
   where \( Y\widehat\times_{\Spec k}\st \) means the formal completion 
   of \( Y\times_{\Spec k}\st \) along the subscheme 
   \( Y\times _{\Spec k} \{0\} \).
\end{prop}

\begin{say}
  A morphism \( \Phi:X\to Z \) of varieties over \( k \) induces  
   a canonical morphism \( \Phi_{m}:X_{m}\to Z_{m} \) 
    \( (m\in \bN\cup \{\infty\}) \).
  Some properties of \( \Phi \) are inherited by \( \Phi_{m} \); 
  for example, if \( \Phi \) is a closed immersion, an open immersion 
  or \'etale, then \( \Phi_{m} \) is also a closed immersion, an open 
  immersion or \'etale.
  But many properties of \( \Phi \) are not inherited by \( \Phi_{m} 
  \); 
  for example, properness, projectiveness, closedness, and so on.
\end{say}

Next we study the  jet schemes and the arc 
space of a variety which admits an action of a group scheme.

\begin{prop}
  Let \( G \) be a group scheme of finite type over \( k \).
  Then \( G_{m} \) \( (m\in \bN\cup \{\infty\}) \) is again a group 
  scheme over \( k \).
  If \( G \) is irreducible, then \( G_{m} \) is also irreducible.  
\end{prop}

\begin{pf}
   Let \( \mu:G\times G\to G \) be the multiplication of the group, 
   let \( e\in G \) be the unit element of the group and let \( 
   \iota:G\stackrel{\sim}\longrightarrow G \) be the morphism defining the inverse elements.
   Then, \( G_{m} \) becomes a group scheme with 
   \( \mu_{m}:G_{m}\times G_{m}\to G_{m} \) the multiplication of the 
   group,
   where \( \mu_{m} \) is induced on \( (G\times G)_{m}\simeq  G_{m}\times G_{m}\)
    from \( \mu \).
  The  scheme \( \{e\}_{m} \)  is a \( k \)-valued point of \( G_{m} \) and it 
  is the unit element under this multiplication.
  The morphism \( \iota_{m}:G_{m}\stackrel{\sim}\longrightarrow G_{m} \) induced from \( \iota   \)
  gives the inverse elements.
  If \( G \) is irreducible, then it is a non-singular irreducible 
  variety which yields that \( G_{m} \) is also non-singular and irreducible.            
\end{pf}

\begin{prop}
  Let \( G \) be a group scheme of finite type over \( k \) and 
  \( X \)  a variety admitting an action of \( G \).
  Then, for \( m \in \bN\cup \{\infty\} \), \( X_{m} \) admits a 
  canonical action of \( G_{m} \) induced from the action of \( G \) 
  on \( X \).  
\end{prop}

\begin{pf}
  Let \( \psi:G\times X\to X \) be the morphism defining the action 
  of \( G \) on \( X \).
  Then the morphism \( \psi_{m}:G_{m}\times X_{m}\simeq (G\times X)_{m}
  \to X_{m} \) induced 
  from \( \psi \) gives an action of \( G_{m} \) on \( X_{m} \).
\end{pf}

\begin{exmp}
  If \( G \) is an \( n \)-dimensional torus \( T^n\simeq 
  (\bA_{k}^1\setminus \{0\})^n \),
  then \( G_{m}\simeq T^n\times \bA_{k}^{nm} \).
  Let \( x=(x^{(0)}_{1},\ldots,x^{(0)}_{n},x^{(1)}_{1},\ldots,x^{(1)}_{n},
  \ldots,x^{(m)}_{1},\ldots,x^{(m)}_{n}) \) and 
  \( y=(y^{(0)}_{1},\ldots,y^{(0)}_{n},y^{(1)}_{1},\ldots,y^{(1)}_{n},
  \ldots,y^{(m)}_{1},\ldots,y^{(m)}_{n}) \) be two \( k \)-valued 
  points of \( G_{m} \),
  where \( (x^{(0)}_{1},\ldots,x^{(0)}_{n}), (y^{(0)}_{1},\ldots,y^{(0)}_{n})
  \in T^n \).
  Then the multiplication \( x\cdot y \)  of \( x \) and \( y \) is 
  \( (x^{(0)}_{1}y^{(0)}_{1},\ldots,x^{(0)}_{n}y^{(0)}_{n},
  \sum_{i+j=1}x^{(i)}_{1}y^{(j)}_{1}, ...,\\
  \sum_{i+j=1}x^{(i)}_{n}y^{(j)}_{n}, \ldots,
  \sum_{i+j=m}x^{(i)}_{1}y^{(j)}_{1}, \ldots,
  \sum_{i+j=m}x^{(i)}_{n}y^{(j)}_{n}) \).
  The unit element of \( G_{m} \) is   
 \[  (\overset{\text{\(  n \) times}}{\overbrace{1,\ldots,1}} ,0,\ldots,0).    \]
 
\end{exmp}

\begin{exmp}
  Let \( X \) be a toric variety with the torus \( T \).
  Then \( T_{m} \) acts on \( X_{m} \) for every \( m\in \bN\cup 
  \{\infty\} \).
\end{exmp}

\begin{say}
  As the \( m \)-jet scheme \( X_{m} \) of a variety \( X \) is of 
  finite type over \( k \), a point of \( X_{m} \) is closed if and 
  only if it is a \( k \)-valued point. 
  But \( X_{\infty} \) is not of finite type and the equivalence above 
  does not hold. 
   First we will see the affirmative case under a condition on \( k \). 
\end{say}

\begin{prop}
  Assume that the base field \( k \) is uncountable.
  Then, for every variety \( X \), a  point of \( X_{\infty} \) 
  is closed 
  if and only if the point is a \( k \)-valued point.
\end{prop}

\begin{pf}
  As the problem is local, we may assume that \( X \) is affine.
  Therefore we have only to prove the assertion for the case   
  \( X_{\infty}=\Spec R \), 
  \( R= k[x_{1},x_{2}, \ldots, x_{n}, \ldots] \), 
  where the variables \( x_{1},x_{2},\ldots,x_{n},\ldots \) are 
  countably infinite.
  For the assertion of the proposition, 
  it is sufficient to prove that  every prime ideal \( I\subset 
  k[x_{1},x_{2},\ldots,x_{n},\ldots] \) is contained in 
  a maximal ideal \( (x_{1}-a_{1}, x_{2}-a_{2}, 
  \ldots,x_{n}-a_{n},\ldots) \), \( a_{1},a_{2},\ldots,a_{n},\ldots\in 
  k \).
  For every \( n \), let \( R_{n} \) be a subring  \( k[x_{1},\ldots,x_{n}]\) of \( 
  R \) and \( I_{n} \) be the intersection \( I\cap R_{n} \). 
  For \( m< n \) the inclusion \( R_{m}\hookrightarrow R_{n } \)   induces the 
  projection \( \Spec R_{n }\to \Spec R_{m} \) which 
  induces a dominant map 
  \( \psi_{n,m}: Z(I_{n})\to Z(I_{m}) \), 
   \( (a_{1},\ldots,a_{m},\ldots,a_{n })\mapsto (a_{1},\ldots a_{m}) \),
   where \( Z(I_{n}) \) is the set of closed points of the closed 
   subscheme defined by \( I_{n} \).
  Fix \( r\geq 1 \). 
  Since \( Z(I_{r})\neq \emptyset \)
  for every \( n>r \), \( \image \psi_{n,r} \) is a non-empty 
  constructible set  and
  \[\image \psi_{r+1, r}\supset \image \psi_{r+2,r}\supset\ldots  \]
  is a non-increasing sequence.
  As \( k \) is uncountable, the intersection \( \bigcap_{n>r} \image \psi_{n,r}\)
  is non-empty by \cite[Proposition 6.5]{batyrev}.
  Take a point \( p_{r} \) from this set.
  In \( Z(I_{r+1}) \),        

 \[ \psi_{r+1,r}^{-1}(p_{r})\cap \image \psi_{r+2,r+1}\supset
    \psi_{r+1,r}^{-1}(p_{r}) \cap \image \psi_{r+3,r+1}\supset\ldots
    \]
  is   a non-increasing sequence of non-empty constructible sets.
  Therefore, we can take a point \( p_{r+1}\in 
  \psi_{r+1,r}^{-1}(p_{r})\cap \left( \bigcap_{n>r+1} \image 
  \psi_{n,r+1} \right) \).
  In the same way, we have points \( p_{r+2}\in Z(I_{r+2})\), \(p_{r+3} 
  \in Z(I_{r+3}),\ldots \) such that \( \psi_{n+1,n}(p_{n+1})=p_{n} 
  \in Z(I_{n})\) 
  for \( n\geq r \).     
Therefore, there is a sequence
  \( a_{1},a_{2},\ldots,a_{n},\ldots \in k \) 
  such that \( p_{n}=(a_{1},a_{2},\ldots,a_{n}) \).
Hence,  
  \( I_{n}\subset (x_{1}-a_{1}, x_{2}-a_{2}, 
  \ldots,x_{n}-a_{n}) \) for every \( n \).
  Then, it follows \( I=\varinjlim I_{n}\subset (x_{1}-a_{1}, x_{2}-a_{2}, 
  \ldots,x_{n}-a_{n},\ldots)  \).  
\end{pf}

\noindent
In the proposition above, 
the condition on \( k \) is essential.
In fact, we obtain the following:

\begin{prop}[Watanabe, Yoshida]
\label{countable}
  Let \( k \) be a countable  field.
  Then there is a closed point which is not a \( k \)-valued point 
  in \(\Spec k[x_{1}, x_{2},  \ldots,  x_{n},  \ldots] \).
\end{prop}

\begin{pf}
  Let \( y \) be a transcendental element over \( k \).
  As \( k \) is countable, the extension field \( k(y) \) is a countably generated \( k 
  \)-algebra.
  Therefore there exists a surjective homomorphism 
  \( k[x_{1}.x_{2},\ldots,x_{n},\ldots]\to k(y) \).
  The kernel of this homomorphism is a maximal ideal which does not 
  give a \( k \)-valued  point.
\end{pf}

As we assume that the base field is an arbitrary algebraically closed 
field, a closed point of an arc space  is not 
necessarily a \( k \)-valued point.
In spite of such a difficulty 
we can see  the structure of the arc space for a toric variety.

\section{Basic properties of the arc space of a toric 
variety}

\begin{say}
  We use the notation and terminology of \cite{fulton}.
  Let  $M$  be the free abelian group  ${\bZ}^n$ $(n\geq 1)$
  and  $N$   its dual $\Hom_{\bZ}(M, {\bZ})$.
  We denote  $M\otimes _{\bZ}{\bR}$  and $N\otimes_{\bZ}{\bR}$  by
  ${M_{\bR}}$  and  $N_{\bR}$, respectively.
  The canonical pairing \( \langle\ , \ \rangle:
  N\times M \to \bZ \)  extends to  
  \( \langle\ \ , \ \rangle:
  N_{\bR}\times M_{\bR} \to \bR \).
  For a linear subspace \( W\subset N_{\bR} \),
  the induced pairing 
  \( (N_{\bR}/W)\times W^{\perp}\to \bR \) 
  is also denoted by \( \langle\ \ , \ \ \rangle \).
  Here,  for \( v\in N_{\bR} \), \( u\in W^{\perp} \) 
  we have that \( \langle v, u \rangle= \langle \rho(v), u \rangle 
  \), 
  where \( \rho: N_{\bR}\to N_{\bR}/W \) is the projection.
  
  For a finite fan ${\D}$ in ${N}$, the corresponding 
   toric variety is denoted by  $T_{N}({\D})$.
  If \( \D \) is the fan consisting of all faces of a cone \( \sigma \), 
  then \( T_{N}(\D) \) is affine and sometimes denoted by \( 
  T_{N}(\sigma) \). 
  
  For a cone \( \tau\in \D \) 
  we denote by $U_{\tau}$  the invariant affine open subset which contains
  ${orb\ 
  {\tau}}$ as the unique closed orbit.
   The open set \( U_{\tau} \) is isomorphic to \( T_{N}(\tau ) \). 
   
   We can write $k[M]$  as  $k[x^{u}]_{u\in M}$,
  where we use the shorthand $x^{u}=x_1^{u_1}x_2^{u_2}\cdots x_n^{u_n}$
  for ${u}=(u_1,\ldots ,u_n)\in M$.
  The torus \( \Spec k[M] \) is denoted by \( T \). 
  We also write \( T \) for the open orbit of the toric variety.
  
\end{say}
   
\begin{prop}
\label{surj}
  Let \( X \) be a toric variety over \( 
  k \) and \( f:Y\to X \) an equivariant resolution of the singularities.
  Then, the induced morphism \( f_{\infty}:Y_{\infty}\to X_{\infty} 
  \) is surjective in a strong sense; i.e. for every extension field 
  \( K\supset k \) the corresponding morphism 
  \( Y_{\infty}(K) \to X_{\infty}(K) \) is surjective.
\end{prop}   
 
\begin{pf}
   Let \(\alpha: \sT \to X \) be an arc of \( X \), then the generic point \( 
  \eta\in \sT \) is mapped to \( orb \tau \) for some cone \( \tau 
   \) in the defining fan of \( X \). 
  As \( f \) is equivariant, \( f^{-1}(orb \tau) \) contains a 
  subscheme isomorphic to \( orb \tau \times T^s \), where \( T^s \) 
  is the torus of dimension \( 0\leq s <n \).
  Hence the restriction \( \Spec K((t))\to X \) of \( \alpha \) 
  can be lifted to \( Y \).
  Therefore, by the properness of \( f \), 
  \( \alpha \) can be  lifted to \( Y \).     
\end{pf}

\noindent 
The irreducibility of the arc space of a variety is known for a 
base field of characteristic zero (\cite{kln}).
In  the positive characteristic case,  
\cite[Example 2.13]{i-k} gives an example of non-irreducible arc space.
But for a toric variety, 
the characteristic is not a problem. 

\begin{cor}
\label{irred}
  The arc space of a toric variety \( X \) is irreducible.
\end{cor}

\begin{cor}
  Since the arc space \( X_{\infty} \) of a toric variety contains 
  \( T_{\infty} \) as an open orbit, \( X_{\infty} \) is an almost 
  homogeneous space by the action of \( T_{\infty} \).
\end{cor}

\begin{pf}
  This follows immediately from the irreducibility of \( Y_{\infty} 
  \) and Proposition \ref{surj}.
\end{pf}

\begin{say}
 \label{good}
 An irreducible component of the fiber \( \pi^{-1}(\sing X) \) of 
 the singular locus \( \sing X \subset X\) 
 is called a {\it good component } if it contains an arc \( \alpha \) 
 such that \( \alpha(\eta) \) is in the non-singular locus (\cite{i-k}).
 If the characteristic of the base field is zero, then every 
 component of \( \pi^{-1}(\sing 
 X) \) is a good component, while 
 there is a non-good component for a positive characteristic 
 case (\cite[Example 2.13]{i-k}).
 The following shows that the characteristic does not affect  on this problem
 for a 
 toric variety.
\end{say}

\begin{prop}
  For a toric variety \( X \), every component of \( \pi^{-1}(\sing X) \)
  is a good component.
\end{prop}

\begin{pf}
   Let \( C \) be a non-good component of \( \pi^{-1}(\sing X) \).
   Let \( f:Y\to X \) be an equivariant resolution of the singularities and 
   \( E_{i} \) \( (i=1,2,\ldots,r) \) be the irreducible components 
   of \( f^{-1}(\sing X) \).
   Then, \( \pi_{Y}^{-1}(E_{i}) \)'s are the  irreducible components 
   of \( f^{-1}_{\infty}( \pi^{-1}(\sing X) ) \), where \( 
   \pi_{Y}:Y_{\infty}\to Y \) is the canonical projection.
   By the surjectivity of \( f_{\infty} \) proved in Proposition 
   \ref{surj}, there is a component  \( \pi_{Y}^{-1}(E_{i}) \) mapped 
   to \( C \).
  However, \( \pi_{Y}^{-1}(E_{i}) \) contains an arc whose image of 
  the generic point corresponds to a point in the non-singular locus on \( X \), 
  which is a contradiction.    
\end{pf}

\noindent
  Now we are going to make a stratification of the arc space
   of a toric variety according to the fan.
  From now on we assume that a toric variety \( X \) is defined by a 
  fan \( \D \).
  Let \( X(\tau)\subset X \) be the closure \( \overline{orb \tau} \)
   for the cone \( \tau\in \D \).
  Then \( X(\tau) \) is again a toric variety.

\begin{defn}
  Let \( X \) be a toric variety corresponding to a fan \( \D \).
 We define \( X_{\infty}(\tau) \) as follows:
  \[  X_{\infty}(\tau)=\{\alpha\in X_{\infty} \mid \alpha:\sT \to X 
  \text{\ factors\  
  through\ }
    X(\tau) \]
   \[ \text{but\ does\ not\ factor\ through\ }
   X(\gamma)\ \text{for\ } \gamma \not< \tau\}  \]
\end{defn}  

\begin{rem}
\begin{enumerate}
\label{remark}
\item[(i)]
By definition, we have:
\[ X_{\infty}(\tau)=\{\alpha\in X_{\infty}\mid \alpha(\eta)\in orb 
\tau\} .\]
In particular, 
\[ X_{\infty}(0)=\{\alpha\in X_{\infty}\mid \alpha(\eta)\in 
T \} .\]
\item[(ii)]
 \( X_{\infty}(\tau)=X(\tau)_{\infty}(0) \), where \( 0 \) is the cone 
 consisting of the origin.
\item[(iii)]
\( X_{\infty}  \) is the disjoint union:
\[ X_{\infty}=\bigsqcup_{\tau\in \D}X_{\infty}(\tau). \] 
\end{enumerate}
\end{rem}
\begin{prop}
\label{closure}
  Let \( X \) be a toric variety defined by a fan \( \D \),
  \( T \) the torus acting on \( X \) and \( \tau \)  a cone in \( 
  \D \).
  Then, the subset \( X_{\infty}(\tau) \) is a   locally closed subset which is invariant
  under the action of \( T_{\infty} \).
\end{prop}

\begin{pf}
    As \( X(\gamma) \) is  closed in \( X \) for every cone \( \gamma\in \D \),
   \( 
  X(\gamma)_{\infty} \) is considered as a closed subscheme of \( X_{\infty} \).
  By definition  
  \[(\ref{closure}.1)\ \ \ \ \ \ \ \ \ \ \ \ \ \ \ \  X_{\infty}(\tau)=X(\tau)_{\infty}\setminus \left( \bigcup_{\gamma \not< \tau} 
  X(\gamma)_{\infty}\right)\ \ \ \ \ \ \ \ \ \ \ \ \ \  \]
   as subsets in \( X_{\infty} \), 
   which shows that  \( X_{\infty}(\tau) \) 
  is locally closed.

  As \( X(\gamma) \) is invariant under the action of  \( T \) 
   for every \( \gamma\in \D \),
  \( X(\gamma)_{\infty} \) is invariant under the action of \( T_{\infty} \).
  The description of \( X_{\infty}(\tau) \) as above gives the 
  assertion of the invariance.
\end{pf}
 
\begin{prop}
  Let \( X \) be a toric variety defined by a fan \( \D \)
  and \( \tau \), \( \gamma \) be cones in \( \D \).
  Then, \( \gamma< \tau \) if and only if \( 
  \overline{X_{\infty}(\gamma)}\supset X_{\infty}(\tau) \).
\end{prop}

\begin{pf}
  First note that \( X(\gamma)_{\infty} \) and \( X(\tau)_{\infty} \) 
  are irreducible (Corollary \ref{irred}) and closed in \( X_{\infty} \).
  Then, the description (\ref{closure}.1) gives that 
  \( \overline{X_{\infty}(\gamma)}=X(\gamma)_{\infty} \) and
  \( \overline{X_{\infty}(\tau)}=X(\tau)_{\infty}  \).
  Therefore, the relation \( \overline{X_{\infty}(\gamma)}
  \supset X_{\infty}(\tau) \) holds if and only if 
   \( X(\gamma)_{\infty}\supset X(\tau)_{\infty} \) holds,
   which is equivalent to \( X(\gamma)\supset X(\tau) \).
  It is well known that the last relation is equivalent to \( \gamma<\tau \).    
\end{pf}


\section{Orbits on the arc space of a toric variety}

\noindent
In this section we associate each \( T_{\infty} \)-orbit on \( 
X_{\infty} \) to a lattice point,
 and describe the dominant relation of two orbits in terms of the 
 corresponding lattice points.

\begin{thm}
\label{orbit}
  Let \( X \) be a toric variety defined by a fan \( \D \). Then,
\begin{enumerate}
  \item[(i)]
  there is a surjective canonical map
  \[ \Psi:X_{\infty}(0)\to |\D |\cap N,\ \ \alpha \mapsto v_{\alpha}, \]
  \item[(ii)]
    for every \( v\in |\D |\cap N  \) there exists a \( k \)-valued 
    point \( \alpha\in X_{\infty}(0) \) such that 
    \[ \Psi^{-1}(v)=T_{\infty}\cdot \alpha, \]
    where \( T_{\infty}\cdot \alpha \) is the orbit of \( \alpha \) 
    by the action of \( T_{\infty} \), and
  \item[(iii)]
      for \( v\in  |\D |\cap N   \), \( \Psi^{-1}(v) \) is a locally closed 
      subset of \( X_{\infty} \).     
\end{enumerate}  
\end{thm}
  
\begin{pf}
  For a \( K \)-valued point \( \alpha\in X_{\infty}(0) \), take a  cone \( 
  \sigma \in \D \) such that \( \alpha(0)\in U_{\sigma} \).
  Then  \( \alpha \) is an arc of \( U_{\sigma} \)  with \( 
  \alpha(\eta)\in T \), therefore we have a 
  commutative diagram:
  \[ \begin{array}{c c c}
  k[\sigma^{\vee}\cap M]&\stackrel{\alpha^*}\longrightarrow & K[[t]]\\
  \bigcap & & \bigcap\\
  k[M] & \stackrel{\alpha^*}\longrightarrow & K((t))\\
  \end{array} \] 
  Let \( v_{\alpha}:M\to \bZ \) be a map defined by \( u \mapsto
  \ord \alpha^*(x^u) \).
  Then \( v_{\alpha} \) is a group homomorphism, therefore \( 
  v_{\alpha}\in N \) with the pairing \( \langle v_{\alpha}, u\rangle = 
  \ord \alpha^*(x^u) \).
  For \( u\in \sigma^{\vee}\cap M \), it follows \( \langle v_{\alpha}, u\rangle = 
  \ord \alpha^*(x^u)\geq 0 \), which implies that \( v_{\alpha}\in 
  \sigma \).
  Now we obtain a map \( \Psi :X_{\infty}(0)\to |\D |\cap N \), \( 
  \alpha\mapsto v_{\alpha} \).
  To show the surjectivity, take a point \( v\in |\D |\cap N \).
  Let \( \sigma \) be a  cone containing \( v \).
  Let \( \alpha^*:k[M]\to k((t)) \) be a \( k \)-algebra homomorphism 
  defined by \( \alpha^*(x^u)=t^{\langle v, u\rangle} \) for \( u\in M \).
  Then, \( \alpha^*(k [\sigma^{\vee}\cap M])\subset k[[t]]\), since 
  \( \langle v, u \rangle \geq 0 \) for \( u\in \sigma^{\vee} \).
  Hence, \( \alpha^* \) gives a \( k \)-valued point \( \alpha \) in \( 
   X_{\infty}(0) \).
   
   For (ii), we prove the equality \( \Psi^{-1}(v)=T_{\infty}\cdot \alpha 
   \) for a \( k \)-valued point \( \alpha\in \Psi^{-1}(v) \).
  For a \( k \)-valued point \( \alpha\in X_{\infty}(0) \), 
  take a cone \( \sigma \) such that \( \alpha\in (U_{\sigma})_{\infty} \).
  Then \( \alpha \) corresponds to a ring homomorphism 
  \( \alpha^*: k[\sigma^{\vee}\cap M]\to k[[t]] \).
  On the other hand, a \( K \)-valued point \( \gamma\in T_{\infty} 
  \) 
  corresponds to a ring homomorphism 
  \( \gamma^* : k[M] \to K[[t]] \).
  This homomorphism is equivalent to a ring homomorphism 
  \( \gamma^*: k[\sigma^{\vee}\cap M]\to K[[t]]\) 
  such that the order of \( \gamma^*(x^u) \) is zero for every \( u\in 
  \sigma^{\vee}\cap M \),
  because \( \sigma^{\vee}\cap M \) generates \( M \).
  Then, \( \gamma \cdot \alpha \) corresponds to the homomorphism 
  \( k[\sigma^{\vee}\cap M]\to K[[t]] \) which maps \( x^u \) to 
  \( \gamma^*(x^u)\alpha^*(x^u) \).
  
  Now let \( \alpha\in (U_{\sigma})_{\infty} \) be the arc 
  corresponding to \( v \) which was constructed in (i).
  If \( \beta\in T_{\infty}\cdot \alpha \), then there exists 
  a \( K \)-valued point \( \gamma\in T_{\infty} \) such that 
  \( \beta=\gamma\cdot \alpha \).
  Then, by the above remark, it follows that \( \beta \in 
  (U_{\sigma})_{\infty} \) and \( \beta \) corresponds to 
  \( \beta^*: k[\sigma^{\vee}\cap M]\to K[[t]]\) which maps 
  \( x^u \) to \( \gamma^*(x^u)t^{\langle v, u \rangle} \)
  whose order is \( {\langle v, u \rangle} \).
  Therefore \( \beta\in \Psi^{-1}(v) \).
  Conversely, suppose that \( \beta \in \Psi^{-1}(v)\) and let 
  \( \sigma \) be a cone such that \( \beta \in 
  (U_{\sigma})_{\infty} \).
  Then we can define \( \gamma\in T_{\infty} \) by 
  \( \gamma^*:  k[\sigma^{\vee}\cap M]\to K[[t]]\),
  \( \gamma^*(x^u)=t^{-\langle v, u\rangle}\beta^*(x^u) \).
  For this \( \gamma \) we have  that \( \gamma\cdot \alpha= \beta \).

  For the assertion (iii), take a cone \( \sigma\in \D \) such that 
  \( T_{\infty}\cdot \alpha\subset (U_{\sigma})_{\infty} \).
  It is sufficient to prove that \( T_{\infty}\cdot \alpha \) is 
  locally closed in \( (U_{\sigma})_{\infty}\cap X_{\infty}(0) \).
  Denote \( (U_{\sigma})_{\infty}  \) by \( \Spec A \).
  Let \( \Lambda: k[\sigma^{\vee}\cap M]\to A[[t]] \) be the ring 
  homomorphism induced from the universal family of arcs on 
  \( (U_{\sigma})_{\infty}  \) (see Proposition \ref{ft}).
  Let \( \Lambda(x^{u_{j}})=\sum_{i\geq 0}a_{j,i}t^i \) for   
  generators \( u_{j} \) \( (j=1,\ldots,r) \) of the semigroup  
   \( \sigma^{\vee}\cap M \).
   Then \[  T_{\infty}\cdot \alpha=\Psi^{-1}(v)=
   \{\beta\in  (U_{\sigma})_{\infty}\cap X_{\infty}(0) \mid
   a_{j,i}(\beta)=0 \operatorname{\ for}\ i<\langle v, u_{j}\rangle,\]
  \[\ \ \ \ \ \ \ \ \ \ \ \ \ \ \ \ \ \ \ \ \ 
   a_{j,i}(\beta)\neq 0 \operatorname{\ for}\ i= \langle v, u_{j}\rangle, 
   \ j=1,\ldots,r \}.  \] 
  Hence,  \(  T_{\infty}\cdot \alpha \) is locally closed in 
   \( (U_{\sigma})_{\infty}\cap X_{\infty}(0). \)   
\end{pf}

\begin{say}
   For \( \tau\in \D \), \( X(\tau) \) is   a toric variety \( 
   T_{N_{\tau}}(\D_{\tau}) \), 
    where \( \D_{\tau} \) consists of the 
   cones \( \overline{\sigma}\subset  N_{\bR}/\tau \bR \) 
   which are the images of the cones \( \sigma \in \D \) such that \( 
   \tau <\sigma \) and \( N_{\tau} \)  is the image 
   of \( N \) in \(  N_{\bR}/\tau \bR  \). 
   The affine open subset \( U_{\overline{\sigma} }\subset X(\tau) \) is
   \( \Spec k[\tau^{\perp}\cap \sigma^{\vee}\cap M] \).
   
   Since \( X_{\infty}(\tau)=X(\tau)_{\infty}(0) \) as is seen in Remark 
   \ref{remark},
   we obtain the following from Theorem \ref{orbit}:  
\end{say}

\begin{cor}
\label{tau}
  Let \( X \) be a toric variety defined by a fan \( \D \) and \( 
  \tau \in \D \). Then,
\begin{enumerate}
  \item[(i)]
  there is a surjective canonical map
  \[ \Psi:X_{\infty}(\tau)\to |\D_{\tau} |\cap N_{\tau},\ \ \alpha \mapsto v_{\alpha}, \]
  \item[(ii)]
     for every \( v\in |\D_{\tau} |\cap N_{\tau}  \) there exists a \( k \)-valued 
    point \( \alpha\in X_{\infty}({\tau}) \) such that 
    \[ \Psi^{-1}(v)=T_{\infty}\cdot \alpha, \]
    where \( T_{\infty}\cdot \alpha \) is the orbit of \( \alpha \) 
    by the action of \( T_{\infty} \), and
  \item[(iii)]
      for \( v\in  |\D_{\tau} |\cap N_{\tau}   \), \( \Psi^{-1}(v) \) is a locally closed 
      subset of \( X_{\infty} \).     
\end{enumerate}  
\end{cor} 

\begin{cor}
\label{union}
\begin{enumerate}
\item[(i)]
 \[   X_{\infty }=\bigcup_{\alpha :\text{\( k \)-valued\ point\ 
  of\ }X_{\infty}}T_{\infty}\cdot \alpha . \]
\item[(ii)]
  For every cone \( \tau \), there is a bijection:
  \[ \{T_{\infty}\cdot \alpha \mid \alpha \text{\ is\ a\ \( k 
  \)-valued \ point\ }\in X_{\infty}(\tau)\}\simeq |\D_{\tau}|\cap 
  N_{\tau}. \]

\end{enumerate}
\end{cor}

\begin{defn}
  As an orbit of a \( k \)-valued point \( \alpha \) in \( X_{\infty}(\tau) \)
  is determined by the lattice point  \( v=v_{\alpha} \in |\D_{\tau}| \), 
  we sometimes denote the orbit \( T_{\infty}\cdot \alpha \) by 
  \( T_{\infty}(v) \).  
\end{defn}

\begin{defn}
  Let \( \sigma \) be a cone in \( N \) and \( v, v' \) two 
  points in \( \sigma \).
  We denote \( v\leq_{\sigma} v' \) if \( v'\in v+\sigma \). 
  It is clear that \( \leq_{\sigma} \) is an order in \( \sigma \).
\end{defn}
  
  Now we are going to study the dominant relation between orbits.
  
\begin{prop}
\label{domination}
  Let \( X \) be a toric variety defined by a fan \( \D \).
  Let \( \alpha\in X_{\infty}(\tau) \) and \( \beta\in 
  X_{\infty}(\gamma) \) be \( k \)-valued points for \( \tau, \gamma 
  \in \D \).
  If \( \overline{T_{\infty}\cdot \alpha}\supset T_{\infty}\cdot \beta  \),
  then \( \tau< \gamma \) and  there exists a cone \( \sigma\in \D \) 
  containing \( \tau \) and \( \gamma \)
  such that \( \alpha, \beta\in (U_{\sigma})_{\infty} \).
\end{prop}  
  
\begin{pf}
  By the condition of the proposition, it follows that \( \beta \in 
  \overline{X_{\infty}(\tau)}=X(\tau)_{\infty} \).
  As \( \beta(\eta)\in orb \gamma \), we have \( orb \gamma\subset 
  X(\tau) \), which implies \(  \tau < \gamma \). 
  To see the second assertion, take a cone \( \sigma\in \D \) such 
  that  \( \beta \in (U_{\sigma})_{\infty} \).
  Then \( 
  \beta(\eta)\in orb(\gamma) \)  implies \( \gamma<\sigma \).
 Since \( (U_{\sigma})_{\infty} \) is an open subset of \( X_{\infty} 
 \) containing \( \beta \), there is an arc \( \alpha' \in T_{\infty}\cdot \alpha \cap 
  (U_{\sigma})_{\infty}  \).
  As \( (U_{\sigma})_{\infty} \) is \( T_{\infty} \)-invariant, 
  it contains both \( T_{\infty}\cdot \alpha \) and \( 
  T_{\infty}\cdot \beta \).     
\end{pf}  

 Hence, in order to interpret the condition of the domination 
\( \overline{T_{\infty}\cdot \alpha}\supset T_{\infty}\cdot \beta \) in 
terms of the corresponding lattice points,  we may assume that \( X \) 
is an affine toric variety.
  If \( X \) is an affine toric variety defined by a cone \( \sigma 
  \) and \( T_{\infty}(v)\subset X_{\infty}(\tau) \) for a face \( \tau<
  \sigma \), then 
  \( v\in \overline{\sigma} \cap N_{\tau} \) by Corollary \ref{tau}, 
  where \( \overline{\sigma}  \) is the image of \( \sigma
  \subset N_{\bR}\) by the projection \( N_{\bR}\to N_{\bR}/\tau\bR \).

\begin{prop}
\label{dominant}
   Let \( X \) be an affine toric variety defined by a cone \( \sigma 
   \) in \( N \). 
  Then,   two orbits \( T_{\infty}(v)  \) and 
  \( T_{\infty}(v') \) in \( X_{\infty}(0) \) satisfy 
  \( \overline{T_{\infty}(v)}\supset T_{\infty}(v') \)    
  if and only if \( v\leq_{\sigma} v' \). 
\end{prop}

\begin{pf}
  Assume \( \overline{T_{\infty}(v)}\supset T_{\infty}(v') \).
  If \( \langle v, u\rangle > \langle v', u\rangle \) for some 
  \( u\in \sigma^{\vee}\cap M \),
  then 
  \[ T_{\infty}(v)\subset \overline{T_{\infty}(v)}
  \cap \{ \alpha \in X_{\infty}(0)\mid \ord \alpha^*(x^u)\geq
  \langle v', u\rangle +1\}, \]
  where the right hand side is a proper closed subset of 
  \( \overline{T_{\infty}(v)} \).
  This is a contradiction.
  Hence, \(  v\leq_{\sigma} v'  \).

  Next, assume that \(  v\leq_{\sigma} v'  \) for \( v, v'\in 
  \sigma\cap N \).
  To prove the converse, we divide the proof into two steps.
  
  \noindent
{\bf  Step 1.} The case \( X \) is non-singular.

\noindent
  Let \( e_{1},e_{2},\ldots,e_{n}  \) be the basis of \( M \) such 
  that \( e_{1},..,e_{r},e_{r+1}^{\pm 1},\ldots,e_{n}^{\pm 1} \)   
  generate \( \sigma^{\vee} \).
  Define a $k$-algebra homomorphism \( \Phi^*:k[\sigma^{\vee}\cap M]\to k[[\lambda, t]] \) by
  \[ \Phi^*(x^{e_{i}})=t^{\langle v', e_{i}\rangle}+\lambda t^{\langle 
  v, e_{i}\rangle} . \]
Here, note that \( \Phi^*(x^{e_{i}})=1+\lambda \) for \( i\geq r+1 \),
 since \( \langle v, e_{i}\rangle =\langle v', e_{i}\rangle = 0 \) 
 for these \( i \)'s.
  Then, we obtain a morphism \( \Phi:\Spec k[[\lambda]]\to 
  X_{\infty}(0) \) such that \( \Phi(0')\in T_{\infty}(v') \) and 
  \( \Phi(\eta')\in T_{\infty}(v) \),
  where \( 0' \) is the closed point and 
  \( \eta' \) is the generic point of \( \Spec k[[\lambda]] \).
  This implies that \( \overline{T_{\infty}(v)} \) contains a point 
  of \( T_{\infty}(v') \).
  As  \( \overline{T_{\infty}(v)} \) is \( T_{\infty} \)-invariant, 
  it follows that 
\(   \overline{T_{\infty}(v)} \supset T_{\infty}(v') \).

\noindent
{\bf Step 2.} The general case.

\noindent
  Define \( \sigma' \) as the cone generated by \( v \) and \( v'- v \).
  Then, note that \( \sigma'\subset \sigma \) and 
  \[  v\leq_{\sigma'} v' . \]
  Let \( N' \) be the subgroup of \( N \) generated by \( v \),  \( v'- v \)
  and \( v_{1},v_{2},\ldots,v_{s}\in N \),
  where their images \( \overline{v_{1}},\overline{v_{2}},\ldots,\overline{v_{s}}
  \in N/N\cap \sigma'\bR\) are a basis of \( N/N\cap \sigma'\bR \).
  Then, the toric variety \( Z =T_{N'}(\sigma')\)  
   is non-singular and there is a canonical 
  equivariant morphism  \[  \varphi:Z \to X  \] with  the 
  surjective morphism  \( T'\to T \) of the 
  tori.  
  By Step 1, \( \overline{T'_{\infty}(v)}\supset T'_{\infty}(v') \) 
  follows from \( v\leq_{\sigma'}v' \).
   Take \( k \)-valued points \( \alpha, \beta\in Z_{\infty}(0) \)
  such that \( v_{\alpha}=v \), \( v_{\beta}=v' \),
  then
  \( T_{\infty}\cdot\varphi_{\infty}(\alpha)= 
  \varphi_{\infty}(T'_{\infty}\cdot \alpha)\) and 
  \( T_{\infty}\cdot\varphi_{\infty}(\beta)= 
  \varphi_{\infty}(T'_{\infty}\cdot \beta)\).
  Therefore \(  \overline{T_{\infty}(v)}\supset T_{\infty}(v') \) 
  follows from    \( v_{\varphi_{\infty}(\alpha)}=v \), \( 
  v_{\varphi_{\infty}(\beta)}=v' \). 
\end{pf}

  As \( X_{\infty}(\tau)=X(\tau)_{\infty}(0) \) we obtain the 
  following as a corollary of Proposition \ref{dominant}.
  
\begin{cor}
\label{generaldominant}
   Let \( X \) be an affine toric variety defined by a cone \( \sigma 
   \) in \( N \). 
  Then, for a face \( \tau<\sigma \),   two orbits \( T_{\infty}(v)  \) and 
  \( T_{\infty}(v') \) in \( X_{\infty}(\tau) \) satisfy 
  \( \overline{T_{\infty}(v)}\supset T_{\infty}(v') \)    
  if and only if \( v\leq_{\overline{\sigma}} v' \),
  where \( \overline{\sigma}\subset N_{\bR}/\tau \bR \) is the image 
  of \( \sigma \). 
\end{cor}

Next we will see the relation of the orbits in   mutually different 
strata.
To see this we need the following combinatorial lemma:

\begin{lem}
\label{combinatorics}
  Let \( \sigma \) be an \( n \)-dimensional cone in \( N \), where 
   \( n=\dim 
  N_{\bR}  \), and 
  \( \tau \) an \( r \)-dimensional face of \( \sigma \).
  Then, there exist a non-singular \( n \)-dimensional cone \( 
  \sigma_{0} \) in \( N \) and its 
  \( r \)-dimensional face \( \tau_{0} \) such that 
  \( \sigma_{0}\subset \sigma \) and \( \tau_{0}\subset \tau \).
\end{lem}

\begin{pf}
  First, subdivide \( \tau \) into  non-singular cones and 
  take one of \( r \)-dimensional cones as \( \tau_{0} \).
  Take any \( n \)-dimensional cone \( \sigma' \)  in \( N \)
  with the face \( \tau_{0} \)
   inside of \( \sigma \), and then subdivide \( \sigma' \) 
  into a non-singular fan \( \Sigma \) by Danilov's procedure \cite[\S 8]{danilov}.
  As \( \tau_{0} \) is non-singular, it is still in the new fan \( 
  \Sigma \) as a 
  cone. Hence, we can take an \( n \)-dimensional 
  non-singular cone \( \sigma_{0} \) with the face \( \tau_{0} \) in 
  \( \Sigma \).
\end{pf}
 
\begin{prop}
\label{different}
  Let \( X \) be an affine toric variety defined by a cone \( \sigma 
   \) in \( N \). 
  Then, two orbits \( T_{\infty}(v)\subset X_{\infty}(\tau) \), \( 
  T_{\infty}(v')\subset X_{\infty}(\gamma) \) satisfy the relation 
  \( \overline{T_{\infty}(v)}\supset T_{\infty}(v') \)
  if and only if \( \tau < \gamma \) and 
  \( \rho(v)\leq_{\overline{\sigma}} v' \),
  where \( \rho:N_{\bR}/\tau\bR\to N_{\bR}/\gamma\bR \) is the 
  canonical projection and \( \overline{\sigma} \) is the image of 
  \( \sigma \) in  \( N_{\bR}/\gamma\bR  \). 
\end{prop} 

\begin{pf}
  First assume that \( \overline{T_{\infty}(v)}\supset T_{\infty}(v') \).
  Then, we have \( \tau< \gamma \) by Proposition \ref{domination}.
  By the assumption, there is a morphism \( \Phi:\Spec k[[\lambda]]\to 
  X_{\infty}(0) \) such that \( \beta:=\Phi(0')\in T_{\infty}(v') \) and 
  \( \alpha:=\Phi(\eta')\in T_{\infty}(v) \),
  where \( 0' \) is the closed point and 
  \( \eta' \) is the generic point of \( \Spec k[[\lambda]] \).
  As \( \alpha\in X_{\infty}(\tau) \), 
  \( \Phi \) factors through \( X(\tau)_{\infty} \).
  This gives the $k$-algebra homomorphism:
  \[ \Phi^*: k[\tau^{\perp}\cap \sigma^{\vee}\cap M]\to k[[\lambda, t]]. \]
  By using \( \Phi^* \), we obtain 
  \( \ord \alpha^*(x^u)\leq \ord \beta^*(x^u)  \) 
  for \( u\in \tau^{\perp}\cap \sigma^{\vee}\cap M \) 
   in the same way as in the proof of Proposition \ref{dominant}.
  Therefore, for \( u\in \gamma^{\perp}\cap \sigma^{\vee}\cap M 
  \subset \tau^{\perp}\cap \sigma^{\vee}\cap M   \) 
  the inequality \( \langle v, u \rangle =\langle \rho(v), u\rangle \leq \langle v', 
  u\rangle  \) holds. 
  Hence, \(  \rho (v)\leq_{\overline{\sigma}} v'  \).
  
  To prove the converse, assume \(  \rho (v)\leq_{\overline{\sigma}} v'  \).
  Then, it is sufficient to prove that 
  \( \overline{T_{\infty}(v)}\supset T_{\infty}(\rho(v)) \), 
  because \( \overline{ T_{\infty}(\rho(v))}\supset T_{\infty}(v') \) 
  follows from Corollary \ref{generaldominant}.
  To prove  \( \overline{T_{\infty}(v)}\supset T_{\infty}(\rho(v)) \), 
  we may assume that \( \gamma=\sigma \), since \( X_{\infty}(\gamma)=
  X(\gamma)_{\infty}(0) \).
  We also can assume that \( \dim \sigma=n=\dim N_{\bR} \), 
  because if \( \dim \sigma = s<n \), then 
  \( T_{\infty}(v)=T_{\infty}^{n-s}\times T_{\infty}^{s}(v) \), 
  \( T_{\infty}(\rho(v))=T_{\infty}^{n-s}\times T_{\infty}^{s}(\rho'(v)) 
  \),
  where \( \rho': \sigma\bR \to \sigma \bR/\tau\bR \) is the projection
  and \( T^s, T^{n-s} \) are \( s \) and \( (n-s) \)-dimensional 
  tori, 
  respectively.
  So the problem is reduced to proving that 
  \( \overline{T_{\infty}^s(v)}\supset T_{\infty}^s(\rho'(v)) \).
  
  Now, for \( \sigma \) and \( \tau \),
   let \( \sigma_{0} \) and \( \tau_{0} \) be as in Lemma 
  \ref{combinatorics}.
  Let \( e_{1},e_{2},\ldots,e_{n} \) be a basis of \( M \) which  
  generate \( \sigma_{0}^{\vee} \) and \(  e_{1},e_{2},\ldots,e_{r} \)
  \( (r<n) \) generate \( \tau_{0}^{\perp}\cap \sigma_{0}^{\vee} \).
  Let \[ \Lambda^*:k[\sigma_{0}^{\vee}\cap M]\to k[[\lambda]]((t)) \]
  be a $k$-algebra homomorphism defined by
  \[ \Lambda^*(x^{e_{i}})=(\lambda +1)t^{\langle v, e_{i}\rangle} 
  \ \operatorname{for}\ i=1,\ldots,r, \]
  \[ \Lambda^*(x^{e_{i}})=\lambda t^{\langle v, e_{i}\rangle} 
  \ \operatorname{for}\ i=r+1,\ldots,n.  \]
  It is easy to check that \( \Lambda^*(x^u)\in k[[\lambda, t]] \)
  for every \( u\in \sigma^{\vee}\cap M \),
  since \( v\in \sigma \).
  Then, we obtain a morphism \( \Lambda:\Spec[[\lambda]]\to X_{\infty} \).
  For every \( u\in \sigma^{\vee}\cap M \), we have \( \ord 
  \Lambda^*(x^u)=\langle v, u\rangle \), therefore \( 
  \alpha:=\Lambda(\eta')\in T_{\infty}(v)\subset X_{\infty}(0) \),
  where \( 0' \) is the closed point and 
  \( \eta' \) is the generic point of \( \Spec k[[\lambda]] \).
  Since \( \tau^{\perp}=\tau_{0}^{\perp} \), 
   \( \beta:=\Lambda(0'): \spec k[[t]]\to X \) factors 
  through \( X(\tau) \) by the definition of \( \Lambda^* \).
  As the corresponding ring homomorphism \( \beta^* \) is extended to 
  a ring homomorphism \( k[\tau^{\perp}\cap M]\to k((t)) \), 
  it follows that \( \beta(\eta)\in orb \tau \), which implies 
  \( \beta\in X_{\infty}(\tau) \).
  For every \( u\in \tau^{\perp}\cap\sigma^{\vee}\cap M \), we have 
  \( \ord \beta^*(x^u)=\langle v, u\rangle =\langle \rho(v), u\rangle \).
  Therefore \( \beta \in T_{\infty}(\rho(v)) \).
  Hence, it follows that \( \overline{T_{\infty}(v)} \) contains a 
  point of \( T_{\infty}(\rho(v)) \). 
  By the \( T_{\infty} \)-invariance of \( \overline{T_{\infty}(v)} 
  \),
  we obtain \( \overline{T_{\infty}(v)}\supset T_{\infty}(\rho(v))  \).                
\end{pf}
 
 Summing up the propositions  and corollary \ref{domination}, 
 \ref{generaldominant},
 \ref{different}, we obtain the following:
 
\begin{thm}
\label{latticethm}
  Let \( X \) be a toric variety and \( T_{\infty}(v) \) and \( 
  T_{\infty}(v') \) two orbits in \( X_{\infty}(\tau) \) and \( 
  X_{\infty}(\gamma) \), respectively.
  Then the following are equivalent:
  \begin{enumerate}
    \item[(i)]
      \( \overline{T_{\infty}(v)}\supset T_{\infty}(v') \),
    \item[(ii)]
      \( \tau< \gamma \),
     there exists a cone \( \sigma> \gamma \) 
      such that \( {T_{\infty}(v)}, T_{\infty}(v')\subset 
      (U_{\sigma})_{\infty} \),
      and
      \( \rho(v)\leq_{\overline{\sigma}}v' \), 
      where \( \rho:N_{\bR}/\tau\bR \to N_{\bR}/\gamma\bR \) is the 
      projection and \( \overline{\sigma} \) is the image of \( \sigma 
      \) in \(  N_{\bR}/\gamma\bR \).
  \end{enumerate}
\end{thm}

\begin{say}
  By now, the dominant relation of orbits is discussed in terms of 
  the order relation of lattice points.
  This gives a relation between arc families and valuations, 
  which will be discussed in the next section.
  But the dominant relation of orbits can be more simply described in 
  terms of homomorphisms of semigroups.
  
  If \( X \) is an affine toric variety defined by a cone \( \sigma 
  \) and \( T_{\infty}(v)\subset X_{\infty}(\tau) \) for a  face 
  \( \tau < \sigma \),
  then \( v\in \overline{\sigma} \cap N_{\tau} \subset N_{\bR}/\tau\bR \),
  where \( \overline{\sigma} \) is the image of \( \sigma \) in \( 
  N_{\bR}/\tau \bR \).
  Then, \( v \) can be considered as a semigroup homomorphism \( v:\tau^{\perp}\cap 
  \sigma^{\vee}\cap M \to \bZ_{\geq 0} \).
  Here, \( v \) can be extended as a semigroup homomorphism 
  \( v:\sigma^{\vee}\cap M \to \bZ_{\geq 0}\cup \{\infty\} \),
  where we define \( v(u)=\infty \) for every \( u\not\in \tau^{\perp} \).
  
  Conversely, every semigroup homomorphism 
  \( v:\sigma^{\vee}\cap M \to \bZ_{\geq 0}\cup \{\infty\} \) 
  is obtained by such an extension from an element of \(  
  \overline{\sigma} \cap N_{\tau} \subset N_{\bR}/\tau\bR
   \) for some face \( \tau \).
\end{say}

\begin{lem}
  Let \( \sigma  \) be a cone in \( N \) and \( v:\sigma^{\vee}\cap 
  M\to \bZ_{\geq 0}\cup \{\infty\}  \) a homomorphism of semigroups.
  Then, there exists a face \( \tau< \sigma \) such that 
  \( v^{-1}(\bZ_{\geq 0})=\tau^{\perp}\cap \sigma^{\vee}\cap M \).
\end{lem}

\begin{pf}
  Take the minimal face \( \gamma  \) of \( \sigma \) containing \(C 
  =v^{-1}(\bZ_{\geq 0})\).
  Then, \( C \) contains a relative interior point \( u \) of \( \gamma\).  
  We will show that \( C =\gamma \cap M \).
  Assume that there exists a point \( u_{0}\in \gamma \cap M \) such that 
  \( v(u_{0})=\infty \). 
  Then, note that \( u_{0}+\sigma^{\vee}\subset v^{-1}(\infty) \).
  Let \( \sigma^{\vee} \) be generated by \( u_{1},u_{2},\ldots,u_{r} \).
  Then, there is a representation 
  \( u=\sum_{i=1}^r  a_{i}u_{i} \) with \( a_{i}>0 \) for every \( i \)
  and   
   \( u_{0}=\sum_{i=1}^r  b_{i}u_{i} \) with \( b_{i}\geq 0 \) for 
   every \( i \).
   Then, in the equality:
      \[  mu=\sum b_{i}u_{i}+ \sum_{i}(ma_{i}-b_{i})u_{i}, \]
   the second term of the right hand side is in \( \sigma^{\vee} \) 
   for \( m\gg 0 \).
   Hence, \( v(mu)=\infty \), but this contradicts to that \( 
   v(mu)=mv(u)\in \bZ_{\geq 0} \).
  Now, we obtain that  \(  C =\gamma \cap M \) and \( \gamma \) can be written 
  as \( \tau^{\perp}\cap \sigma^{\vee} \) for some \( \tau < \sigma \).      
\end{pf}

By Corollary \ref{tau} and Theorem \ref{latticethm}, we obtain the following 
interpretation:

\begin{thm}
  Let \( X \) be a toric variety defined by a fan \( \D \), 
  then we obtain the following:
  \begin{enumerate}
  \item[(i)]
  There is a bijective map:
  \[ \{T_{\infty}\cdot{\alpha}\mid \alpha:k\text{-valued\ point\ of }X\}
  \isom \bigsqcup_{\sigma}\Hom_{s.g.}(\sigma^{\vee}\cap M, \bZ_{\geq 0}\cup
   \{\infty\}), \]
   where \( \sigma \)   varies the maximal cones 
   in \( \D \).
   Via this map,
   each \( T_{\infty}\cdot \alpha \) can be written as 
  \(T_{\infty}(v) \) for  a suitable element \( 
  v \) of the right hand side.
   
  \item[(ii)]
  We have the relation \( \overline{T_{\infty}(v)}\supset T_{\infty}(v') \) if and 
  only if there is a maximal cone \( \sigma \) in \( \D \) such that 
  \( v, v'\in \Hom_{s.g}(\sigma^{\vee}\cap M, \bZ_{\geq 0}\cup 
  \{\infty\}) \) and \( v\leq v' \), 
  where \( v \leq v'\) means that \(  v(u)\leq v'(u) \) for every \( 
  u\in \sigma^{\vee}\cap M \).
  \end{enumerate}
\end{thm}

\section{Contact loci of an invariant ideal}

\noindent
In this section, 
  we will give the answer to the embedded version of Nash problem for 
  an invariant ideal of a toric variety.
 
\begin{defn}
  Let \( X \) be a variety over an algebraically closed field \( k \) 
  and \( k(X) \) the rational function field of \( X \).
  A {\it divisorial valuation} of \( k(X) \) is a positive integer times discrete 
  valuation \( val_{D} \) associated to a prime divisor \( D \) on some normal variety \( X' 
  \) which is birational to \( X \).
  Note that this definition is wider than the definition of 
  ``divisorial valuation'' in \cite{ELM}.
\end{defn}

\begin{defn}
  Let \( X \) be an affine toric variety defined by a cone \( \sigma \) in \( N \).
  For every point \( v\in \sigma\cap N \) we can associate a valuation 
  \( val_{v} \) on \( k(X) \) as follows:
  
  \noindent 
  Define
  \[ val_{v}(f):= \min _{x^u\in f} \langle v, u \rangle, 
  \operatorname{\ for\ } f\in k[\sigma^{\vee}\cap M] \]
  and extend it on \( k(X) \), the quotient field of  \( k[\sigma^{\vee}\cap M] \).
  This valuation is called a {\it toric valuation}.
  Here \( {x^u\in f} \) means that the coefficient of the monomial \( x^u \) 
  in \( f \) is not zero.
  Note that the  toric valuation defined by a primitive element \( v \) is 
   \( val_{D_{v}} \),
  where   
 \( D_{v} \) is  the irreducible invariant divisor  
 \( \overline{orb(\bR_{\geq 0} v)} \) on some  
  toric variety \( X' \) which is birational to \( X \).
   Since every toric valuation is a positive integer times such a 
   valuation,
    every toric valuation is a divisorial valuation. 
\end{defn}
 
 \begin{say}
   For a variety \( X \) over an algebraically closed field \( k \), 
   let \( \psi_{m}: X_{\infty} \to X_{m} \) \( (m\in \bZ_{\geq 0}) \)
   be the truncation morphism.
   Note that \( \psi_{0}=\pi \).
   Recall that a cylinder \( C \) in \( X_{\infty} \) is a subset of 
   the form \( \psi_{m}^{-1}(S) \), for some \( m \) and some 
   constructible subset \( S\subset X_{m} \).
\end{say}

\begin{exmp}
  Let \( X \) be a toric variety.
  Then an orbit \( T_{\infty}(v) \)  of a \( k \)-valued point in \( 
  X_{\infty}(0) \) is a cylinder.
  Indeed, we may assume that \( X \) is the affine toric variety defined by a 
  cone \( \sigma \).
  The orbit is the subset of \( X_{\infty} \) consisting of arcs 
  \( \alpha  \) whose corresponding homomorphisms 
  \(\alpha^*: k[\sigma^{\vee}\cap M] \to K[[t]] \) satisfy 
  \( \ord \alpha^*(x^{u_{i}})=\langle v, u_{i}\rangle \) for 
  generators \( u_{1},\ldots, u_{s} \) of \( \sigma^{\vee}\cap M \).
  Let \( m\geq \max_{i=1,\ldots,s}\langle v, u_{i}\rangle \) and 
  \( S_{m}\subset  X_{m}  \) the subset consisting of \( m \)-jets 
  \( \gamma \) whose corresponding homomorphisms 
  \( \gamma^*:k [\sigma^{\vee}\cap M] \to \sTm \) satisfy 
  \( \ord \gamma^*(x^{u_{i}})=\langle v, u_{i}\rangle \).
  Then, \( S_{m} \) is a locally closed subset of \( X_{m} \) and 
  \( T_{\infty}(v)=\psi_{m}^{-1}(S_{m}) \).    
\end{exmp}

\begin{say}
\label{valelm}
  Let \( X \) be a non-singular variety over \( \bC \) and \( C \) an 
  irreducible cylinder in \( X_{\infty} \). 
  In \cite{ELM} a valuation \( val_{C} \) corresponding to \( C \) is 
  defined as follows:
  Note first that if \( \alpha\in X_{\infty} \) is a \( \bC \)-valued 
  point,
  and if \( f \) is a rational function on \( X \) defined in a 
  neighborhood of \( \pi(\alpha) \), then \( \ord{\alpha}^*(f) 
  \) is well defined, 
  where \( \alpha^*:\o_{X}\to \bC[[t]] \) is the ring homomorphism 
  corresponding to \( \alpha \).
  If the domain of \( f \) intersects \( \pi(C) \),
  then \( val_{C}(f):=\ord{\alpha}^*(f) \), for general \( \alpha\in C \).
  Then \( val_{C}(f) \) is well defined and can be extended to a 
  valuation of the function field of \( X \).
\end{say}

\begin{prop}[\cite{ELM}]
\label{elm}
   Let \( X \) be a non-singular variety over \( \bC \) and \( C \) an 
  irreducible cylinder in \( X_{\infty} \) which does not dominate \( X \).
  Then \( val_{C} \) is equal with a divisorial valuation. 
\end{prop}

In the proof of Proposition \ref{elm}, the condition that \( X \) is 
non-singular is used.
Therefore, this proposition does not imply 
that  for a cylinder \(C= T_{\infty}(v)\subset X_{\infty}(0) \) on a singular 
toric variety \( X \), the corresponding valuation \( val_{C} \) is 
a divisorial valuation.
  However, the following proposition shows that \( val_{C} \) 
   is a divisorial valuation for \( C= T_{\infty}(v) \).
  
\begin{prop}
\label{val}
   Let \( X \) be a toric variety over an algebraically closed field \( 
   k \) and \( C=T_{\infty}(v)\subset X_{\infty}(0) \), 
   then, \( val_{C}=val_{v} \). 
   In particular \( val_{C} \) is a divisorial valuation.  
\end{prop}  

\begin{pf}
  We may assume that \( X \) is an affine toric variety defined 
  by a cone \( \sigma \).
  It is sufficient to prove that \( val_{C}(f)=val_{v}(f) \) for 
  every element \( f\in k[\sigma^{\vee}\cap M] \).
  Note that \( val_{C}(f)=\ord\alpha^*(f) \) for the generic point 
  \( \alpha\in C \).
  If \( f \) is a monomial \( x^u \) \( (u\in \sigma^{\vee}\cap M )\),
  then by the definition of \( C= T_{\infty}(v) \) we have
  \[  val_{C}(x^u)=\ord \alpha^*(x^u)=\langle v, u\rangle =val_{v}(x^u).  \] 
  For general \( f \), we have
  \[ val_{C}(f)\geq\min_{x^u\in f} val_{C}(x^u)=
  \min_{x^u\in f}\langle v, u\rangle =val_{v}(f).\]
  On the other hand, let \( R_{v} \) is the discrete valuation ring of the 
  divisorial valuation  \( val_{v} \).
  Then there is an indeterminate \( t \) such that the composite 
\[   \beta^*: k[\sigma^{\vee}\cap M]\hookrightarrow R_{v}\hookrightarrow 
  \widehat{R}_{v}\simeq K[[t^e]]\hookrightarrow K[[t]]      
 \]
 satisfies \( \ord \beta^*(f)=val_{v}(f) \) for \( f\in  k[\sigma^{\vee}\cap M] \).
  Here, \( K \) is the residue field of \( R_{v} \) by the maximal 
  ideal and \( e \) is the positive integer such that \( v=ev_{0} \) for a  
  primitive element \( v_{0} \).
  As the arc \( \beta :\sT \to X \) corresponding to \( \beta^* \) is 
  a \( K \)-valued point of \( C \), 
  we obtain the following inequality by the upper semicontinuity  
  \[ val_{C}(f)=\ord \alpha^*(f)\leq \ord \beta^*(f)=val_{v}(f). \]
  Therefore, we obtain \( val_{C}(f)=val_{v}(f) \). 
\end{pf}

Now we recall the definition of the contact locus of an ideal of a variety \( X \).
Let \( X \) be an affine variety over an algebraically closed field \( k \)
with the coordinate ring \( A \) and \( \frak{a} \) an ideal of \( A \).
Then, we define the \( p \)-th contact locus of \( \frak{a} \) by
\[ \Cont^p({\frak{a}})=\{\alpha \in X_{\infty} \mid \min_{f\in 
{\frak{a}}}
\ord \alpha^*(f)=p\}. \]
It is clear that this is a cylinder.
If \( X \) is non-singular then the irreducible components are also
  cylinders.
  Therefore each irreducible component of the contact locus 
  corresponds to a divisorial valuation.
  Now, we can state the embedded version of Nash problem posed in 
  \cite{ELM}.

\begin{problem}
  Which valuations correspond to the irreducible components of 
  \( \Cont^p({\frak{a}}) \)?
\end{problem}  

We consider this problem  for an invariant ideal \( {\frak a} \)
   on a toric variety \( X \).
   We should note that for a singular variety \( X \), an irreducible 
   component of a cylinder is not a cylinder in general, therefore 
   an irreducible component does not necessarily correspond  to a 
   divisorial valuation.
   But in our toric case,  an irreducible component of the contact locus 
   corresponds to a divisorial valuation.
   
\begin{lem}
  Let \( X \) be an affine toric variety and \( {\frak a} \) an 
  invariant ideal on \( X \).
  Then, for every integer \( p>0 \), 
   an orbit \( T_{\infty}(v) \) is either contained in \( \Cont^p({\frak a}) \)
  or  disjoint from \( \Cont^p({\frak a})  \).
\end{lem}  

\begin{pf}
  Take an arc \( \alpha\in T_{\infty}(v) \).
  Then \( \alpha \) belongs to \( \Cont^p({\frak a}) \) if and only 
  if 
  \[  p=\min_{x^u\in {\frak a}}\ord \alpha^*(x^u)=
  \min_{x^u\in {\frak a}}\langle v, u \rangle, \]
  where we define \( \langle v, u \rangle = \infty \) if 
  \( v\in N_{\bR}/\tau\bR \) and \( u\not\in \tau^{\perp} \) 
  for a cone \( \tau \).
  The assertion of the lemma follows immediately from this.
\end{pf}  

By this lemma it follows that \( \Cont^p({\frak a}) \) is a union of 
\( T_{\infty}(v) \)'s.

\begin{lem}
  Let \( X \) be an affine toric variety defined by a cone \( \sigma  \) 
  in \( N \) and \( {\frak a} \) an invariant ideal on \( X \).
  If an orbit \( T_{\infty}(v)\subset \Cont^p({\frak a}) \) is in \( 
  X_{\infty}(\tau) \) for \( \tau\neq 0 \),
  then there is an orbit \( T_{\infty}(\tilde v)\subset X_{\infty}(0)
   \) such that 
    \(  T_{\infty}(\tilde v)\subset\Cont^p({\frak a}) \)
    and 
    \( \overline{T_{\infty}(\tilde 
  v)}\supset T_{\infty}(v) \).
  \end{lem}

\begin{pf}
  Let \( \rho:N_{\bR}\to N_{\bR}/\tau\bR \) be the projection.
  As \( v \) is in the image  \( \rho(\sigma\cap N) \),
  we can take a point \( v_{0}\in \sigma\cap N \) such that 
  \( \rho(v_{0})=v \).
  Then \( \langle v, u\rangle =\langle v_{0}, u\rangle \) for 
  \( u\in \sigma^{\vee}\cap \tau^{\perp} \).
  We can naturally define \( \langle v, u\rangle =\infty  \) for \( u\in 
  \sigma^{\vee}\setminus \tau^{\perp} \). 
  Let \( v_{1}\in \tau\cap N \) be in the relative interior of \( \tau 
  \).
  Then \( \langle mv_{1}, u\rangle > p \) for every 
  \( u\in (\sigma^{\vee}\setminus \tau^{\perp})\cap N \) and an 
  integer 
  \( m>p \).
  Let \( \tilde{v} = v_{0}+mv_{1} \) \( (m>p) \).
  Then, for every \( u\in \tau^{\perp}\cap\sigma^{\vee}\cap M \)
  it follows that 
  \( \langle \tilde v, u\rangle =\langle v_{0}, u\rangle =\langle 
  v,u\rangle \), while
  for every \( u\in (\sigma^{\vee}\setminus \tau^{\perp})\cap M \)
  it follows that \( \langle \tilde v, u\rangle> p \).
  Therefore \[ \min_{x^u \in {\frak a}}\langle \tilde v, u\rangle
  = \min_{x^u\in {\frak a}}\langle  v, u\rangle
  = p. \]
  Hence  \( T_{\infty}(\tilde v)\subset \Cont^p({\frak a}) \).    
  
  On the other hand, 
  \( \rho(\tilde{v})=v \) yields \( \overline{T_{\infty}(\tilde 
  v)}\supset T_{\infty}(v) \), 
  by  Proposition \ref{different}.       
\end{pf}

By these lemmas, we obtain that an irreducible component  of 
\( \Cont^p({\frak a}) \) is the closure of  \( T_{\infty}(v) \) for 
some \( v\in \sigma\cap N \) such that 
\( \min_{x^u\in {\frak a}}\langle v, u\rangle = p \).
Here, by Proposition \ref{dominant} and Proposition \ref{val} , 
  we obtain the answer to 
the embedded version of Nash problem.

\begin{thm}
\label{embNash}
  Let \( {\frak a} \) be an invariant ideal on an affine toric variety \( X \)
  defined by a cone \( \sigma \).
  Then, an irreducible component of \( \Cont^p({\frak a}) \) is 
 the closure of  \( T_{\infty}(v) \) for an element \( v \) minimal
 in \( V({\frak a},p)=\{ v'\in \sigma\cap N \mid 
 \min_{x^u\in {\frak a}}\langle v', u\rangle = p
 \} \) with respect to the order \( \leq_{\sigma} \).
  Therefore the valuations \( \{val_{v}\mid v\in \sigma\cap N \text{\ 
   minimal \  in\ } V({\frak a},p)\} \) correspond bijectively 
    to the irreducible 
  components of \( \Cont^p({\frak a}) \).
\end{thm}

\begin{rem}
  Let \( G({\frak a})\subset M_{\bR} \) be the Newton polytope of
   \( {\frak a} \) as in Figure 1 and \( \D({\frak a}) \) the dual fan of \( G ({\frak a})\).
   The dual fan is the subdivision of \( \sigma \).
   Then, the function \( g(v):=\min_{u\in G({\frak a})}
    \langle v,u\rangle \) \( (v\in \sigma) \) is a strongly convex piecewise linear 
    function with respect to the fan \( \D({\frak a}) \).
    Therefore the subset \( g^{-1}(p)=\{v\in \sigma \mid g(v)=p\}\) is the 
    boundary of some convex polytope as in 
    the Figure 2. 
    The minimal elements of \( V({\frak a},p) \) are on this boundary.
    It is clear that this convex polytope is \( pG({\frak a})^{\circ} 
    \),
    where \( G({\frak a})^{\circ} \) is the polar polytope defined as 
    \( \{v\in \sigma \mid g(v)\geq 1\} \).

    We can see that a lattice 
    point of a compact face of \( g^{-1}(p) \) is always a minimal 
    element of \( V({\frak a},p) \), therefore it gives a valuation 
    corresponding to an irreducible component of \( \Cont^p({\frak a}) \).
    If \( p \) is divisible enough so that every vertex of \( pG({\frak a})^{\circ} 
    \) is in \( N \), 
    then the minimal elements in \( V({\frak a},p) \) coincide with 
    the lattice points on the compact faces 
    of \( g^{-1}(p) \).
    
\end{rem}

\begin{rem}
  The referee kindly informed the following to the author:
  For \( u\in \sigma^{\vee}\cap M \), the log canonical threshold 
  \( {\operatorname{lc}}(X, V({\frak a}), V(x^u)) \) turns out to be the maximal 
  value \( \lambda  \) such that \( x^u\not\in {\cal I}(X, {\frak 
  a}^{\lambda}) \) by \cite{bl}, 
  where \( {\cal I}(X, {\frak 
  a}^{\lambda})  \) is a multiplier ideal for \( {\frak a} \).
  Some multiple of the primitive vector \( v\in \sigma\cap N \) corresponding to a 
  divisor which computes  \( {\operatorname{lc}}(X, V({\frak a}), V(x^u)) \)  
   lies on a compact face of \( g^{-1}(p) \) for some \( p \).
   Conversely, for some multiple of a primitive vector \( v\in \sigma\cap N \) on a 
   compact face of 
   \( g^{-1}(p) \), there exists \( u\in \sigma^{\vee}\cap M  \) such 
   that the divisor corresponding to \( v \) computes the log 
   canonical threshold \( {\operatorname{lc}}(X, V({\frak a}), V(x^u)) \).
\end{rem}

\begin{exmp}
  Let \( X \) be an affine toric variety defined by a cone \( \sigma \).
  Then  the  components in \( \pi^{-1}(\Sing X) \) are 
  \( \overline{T_{\infty}(v)} \)'s,
  where \( v \)'s are the minimal elements in \( \bigcup_{\tau< 
  \sigma:\operatorname{singular}}\tau^o \cap N \) with respect to
  the order \( \leq_{\sigma} \).
  Here, \( \tau^o \) is the relative interior of \( \tau \).
  This is proved as follows: 
  Let  \( {\frak a} \) be the ideal of \( \sing X \), 
  then it is an invariant ideal.
  As \( \pi^{-1}(\sing X)=\bigcup_{p\geq 1}\Cont^p({\frak a}) \),  
  it follows that an irreducible component of \( \pi^{-1}(\sing X) \) 
  is \( \overline{T_{\infty}(v)}  \),
  where \( v \) is minimal among \( v' \)'s such that \( v'\in 
  \sigma\cap N \) and \( \min_{x^u\in {\frak a}}\langle v', u\rangle \geq 1 \)
   by Theorem \ref{embNash}.
  Here, \( \min_{x^u\in {\frak a}}\langle v', u\rangle \geq 1 \) 
  if and only if \( \alpha(0)\in \sing X \)
  for \( \alpha \) with \( v_{\alpha}=v' \), 
  which is equivalent to the fact that \( v'\in \tau^o \) for a 
  singular face \( \tau< \sigma \) by \cite[Proposition 3.9]{i-k}.  
  
\end{exmp}


\makeatletter \renewcommand{\@biblabel}[1]{\hfill#1.}\makeatother

\end{document}